\input amstex
\documentstyle{amsppt}
\magnification1200
\pagewidth{6.5 true in}
\pageheight{9.25 true in}
\NoBlackBoxes

\topmatter
\title What is the best approach to counting primes?
\endtitle
\author Andrew Granville\endauthor

\endtopmatter

\parindent=20pt
\document

\def\phi{\varphi}

\def\Li{\text{\rm Li}}
\def\hbar{\overline{h}}

{\eightpoint As long as people have studied mathematics, they have wanted to know how many primes there are. Getting precise answers is a notoriously difficult problem, and the first suitable technique, due to Riemann, inspired an enormous amount of great mathematics, the techniques and insights permeating many different fields. In this article we will review some of the best techniques for counting primes, centering our discussion around Riemann's seminal paper. We will go on to discuss its limitations, and then recent efforts to replace Riemann's theory with one that is significantly simpler.}\medskip

\noindent {\bf 1. How many primes are there? Predictions} \
You have probably seen a proof that there are infinitely many prime numbers, and were perhaps curious as to roughly how many primes there are up to a given point. With the advent of substantial factorization tables,\footnote{See Appendix 1} it was possible to make predictions supported by lots of data.
On December 24th, 1849, Gauss wrote to his ``most honored friend'', Encke, describing his own attempt to guess at an approximation as to the number of primes up to $x$ (which we will denote throughout by $\pi(x)$).
Gauss describes his work:

\block {\eightpoint
    {\sl First beginning \dots in 1792 or 1793 {\rm [when Gauss was 15 or 16]} \dots I  \dots directed my
attention to the decreasing frequency of prime numbers, to which end I counted
them up in several chiliads {\rm [blocks of 1000 consecutive integers]}
and recorded the results \dots I soon recognized \dots
it is nearly inversely proportional to the logarithm, so that the number of all
prime numbers under a given boundary $x$ were nearly expressed through the
integral
$$ \int_2^x \frac{dt}{\log t} $$
where the logarithm is understood to be the natural logarithm.}  
}
\endblock

\noindent Gauss went on to compare his guess $ \int_2^x \frac{dt}{\log t} $, which we denote by $\Li(x)$, with $\pi(x)$, the actual count of the number of primes up to $x$.
\medskip

{\eightpoint

\centerline{\vbox{\offinterlineskip \halign{\vrule  #&\ \ #\ \hfill
&& \hfill\vrule #&\ \ \hfill # \hfill\ \ \cr \noalign{\hrule} \cr
height5pt&\omit && \omit && \omit & \cr &Under&& $\pi(x) =
\#\{\text{primes} \le x\}=$ &&   $=\Li(x)\pm$  Error & \cr
height5pt&\omit && \omit && \omit & \cr \noalign{\hrule} \cr
height3pt&\omit && \omit && \omit & \cr
&\  $500000$     && 41556 && \, 41606.4 \ - \ \, \ 50.4 &\cr
height3pt&\omit && \omit && \omit & \cr
& $1000000$     && 78501 && \, 79627.5 \ - \ \, 126.5 &\cr
height3pt&\omit && \omit && \omit & \cr
& $1500000$     && 114112 && 114263.1 \ - \ 151.1 &\cr
height3pt&\omit && \omit && \omit & \cr
& $2000000$     && 148883 && 149054.8 \ - \ 171.8 &\cr
height3pt&\omit && \omit && \omit & \cr
& $2500000$     && 183016 && 183245.0 \ - \ 229.0 &\cr
height3pt&\omit && \omit && \omit & \cr
& $3000000$     && 216745 && 216970.6 \ - \ 225.6 &\cr
height3pt&\omit && \omit && \omit & \cr
\noalign{\hrule}}} } \botcaption{Table 1} Primes up
to various points, and a comparison with Gauss's prediction.\endcaption

% \centerline{ (Here, and throughout, we define $\Li(x):=\int_2^x \frac{dt}{\log t}$.) }
\medskip

}

\noindent

In his {\sl Th\'eorie des Nombres}, Legendre proposed
$$
\frac x{\log x-A}
$$
with $A=1.08366$ as a good approximation for $\pi(x)$, in which case the comparative errors are
$$
- 23.3, \ \
+ 42.2, \ \
+ 68.1, \ \
+ 92.8, \ \
+ 159.1, \ \text{and} \
+ 167.6,
$$
respectively. These are smaller than the errors from Gauss's $\Li (x)$, though both seem to be excellent approximations.  Nevertheless Gauss retained faith in his prediction:

\block {\eightpoint
    {\sl
These differences are smaller than those with the integral, though they do
appear to grow more quickly than {\rm [the differences given by the integral]} with
increasing $x$, so that it is possible that they could easily surpass the latter,
if carried out much farther.}}
\endblock

\noindent Today we have data that goes ``much farther'':\medskip

{\eightpoint

\centerline{\vbox{\offinterlineskip \halign{\vrule #&\ \ #\ \hfill
&& \hfill\vrule #&\ \ \hfill # \hfill\ \ \cr \noalign{\hrule} \cr
height5pt&\omit && \omit && \omit && \omit & \cr &$ \ \ x $&& $\pi(x) =
\#\{\text{primes} \le x\}$ &&  Gauss's error term && Legendre's error term&\cr
height5pt&\omit && \omit && \omit && \omit & \cr \noalign{\hrule} \cr
height3pt&\omit && \omit && \omit && \omit & \cr
&  $10^{20}$  && 2220819602560918840 && 222744644 && 2981921009910364 & \cr
height3pt&\omit && \omit && \omit && \omit & \cr
& $10^{21}$ && 21127269486018731928 && 597394254  && 27516571651291205 & \cr
height3pt&\omit && \omit &&\omit && \omit & \cr
& $10^{22}$ && 201467286689315906290 && 1932355208 && 254562416350667927 & \cr
height3pt&\omit && \omit && \omit && \omit & \cr
& $10^{23}$ &&
1925320391606803968923 && 7250186216&& 2360829990934659157 & \cr
height3pt&\omit && \omit&& \omit && \omit & \cr \noalign{\hrule}}} }
\botcaption{Table 2} Comparing the errors in Gauss's and Legendre's predictions.\endcaption
\medskip

}

\noindent It is now obvious that Gauss's prediction is indeed better, that Legendre's error terms  quickly surpass those of Gauss, and keep on growing bigger. Here is some of the most recent data and a comparison to Gauss's guesstimate, $\Li (x)$:\medskip

{\eightpoint

\centerline{\vbox{\offinterlineskip \halign{\vrule #&\ \ #\ \hfill
&& \hfill\vrule #&\ \ \hfill # \hfill\ \ \cr \noalign{\hrule} \cr
height5pt&\omit && \omit && \omit & \cr &$ \ \ x $&& $\pi(x) =
\#\{\text{primes} \le x\}$ &&  Overcount: $[\Li(x)]-\pi(x)$ & \cr
height5pt&\omit && \omit && \omit & \cr \noalign{\hrule} \cr
height3pt&\omit && \omit && \omit & \cr & $10^3$     && 168 && 10 &
\cr height3pt&\omit && \omit && \omit & \cr & $10^4$     && 1229 &&
17 & \cr height3pt&\omit && \omit && \omit & \cr & $10^5$     &&
9592 && 38 & \cr height3pt&\omit && \omit && \omit & \cr & $10^6$ &&
78498 && 130 & \cr height3pt&\omit && \omit && \omit & \cr & $10^7$
&& 664579 && 339 & \cr height3pt&\omit && \omit && \omit & \cr &
$10^8$     && 5761455 && 754 & \cr height3pt&\omit && \omit && \omit
& \cr & $10^9$ && 50847534 && 1701  & \cr height3pt&\omit && \omit
&& \omit & \cr & $10^{10}$  && 455052511 && 3104 & \cr
height3pt&\omit && \omit && \omit & \cr &  $10^{11}$  && 4118054813
&& 11588  & \cr height3pt&\omit && \omit && \omit & \cr &  $10^{12}$
&& 37607912018 && 38263 & \cr height3pt&\omit && \omit && \omit &
\cr &  $10^{13}$ && 346065536839 && 108971 & \cr height3pt&\omit &&
\omit && \omit & \cr &  $10^{14}$  && 3204941750802 && 314890 & \cr
height3pt&\omit && \omit && \omit & \cr &  $10^{15}$  &&
29844570422669 && 1052619 & \cr height3pt&\omit && \omit && \omit &
\cr &  $10^{16}$  && 279238341033925 && 3214632 & \cr
height3pt&\omit && \omit && \omit & \cr &  $10^{17}$  &&
2623557157654233 && 7956589 & \cr height3pt&\omit && \omit && \omit
& \cr &  $10^{18}$  && 24739954287740860 && 21949555 & \cr
height3pt&\omit && \omit && \omit & \cr &  $10^{19}$  &&
234057667276344607 && 99877775 & \cr height3pt&\omit && \omit &&
\omit & \cr &  $10^{20}$  && 2220819602560918840 && 222744644 & \cr
height3pt&\omit && \omit && \omit & \cr & $10^{21}$ &&
21127269486018731928 && 597394254  & \cr height3pt&\omit && \omit &&
\omit & \cr & $10^{22}$ && 201467286689315906290 && 1932355208 & \cr
height3pt&\omit && \omit && \omit & \cr & $10^{23}$ &&
1925320391606803968923 && 7250186216& \cr height3pt&\omit && \omit
&& \omit & \cr \noalign{\hrule}}} }
\botcaption{Table 3} Primes up
to various $x$, and the overcount in Gauss's prediction.\endcaption }
\smallskip

\noindent When looking at this data, compare the widths of the right two columns. The rightmost column is about half the width of the middle column \dots How do we interpret that? Well, the width of a column is given by the number of digits of the integer there, which corresponds to the number's logarithm in base 10. If the log of one number is half that of a second number, then the first number is about the square-root of the first. Thus this data suggests that when we approximate $\pi(x)$, the number of primes up to $x$, by Gauss's guesstimate, $\Li(x)$, the error is around $\sqrt{x}$, which is really tiny in comparison to the actual number of primes. In other words, Gauss's prediction is terrific.

We still believe that Gauss's $\Li(x)$ is {\sl always} that close to $\pi(x)$. Indeed in section 6, we will sketch how the, as yet unproved, {\sl Riemann Hypothesis} implies that
$$|  \pi(x) - \text{Li}(x)|  \leq  x^{1/2} \log x  \tag{RH1}$$
for all $x\geq 3$. This would be an extraordinary thing to prove as there would be many beautiful consequences. For now we will just focus on the much simpler statement that the ratio of $\pi(x) : \Li(x) $ tends to $1$ as $x\to \infty$. Since $\Li(x)$ is well-approximated by $x/\log x$,\footnote{To prove this, try integrating $\Li(x)$ by parts.}  this quest can be more simply stated as
$$
\lim_{x\to \infty} \ \pi(x)\ \big/ \ \frac{x}{\log x} \ \ \text{exists and equals} \ 1.
$$
This is known as the {\sl Prime Number Theorem}, and it took more than a hundred years, and some earth-shaking ideas, to prove it (as we'll outline in sections 4 to 7 of this article).

\subhead 2. Elementary techniques to count the primes\endsubhead It is not easy to find a way to count primes at all accurately. Even proving good upper and lower bounds is challenging.

One effective technique to get an upper bound is to try to use the principle of the sieve of Eratosthenes. This is where we ``construct'' the primes up to $x$, by removing the multiples of all of the primes $\leq \sqrt{x}$.  One starts by removing the multiplies of 2, from a list of all of the integers up to $x$, then the remaining multiples of 3, then the remaining multiples of 5, etc. Hence once we have removed the multiples of the primes $\leq y$ we have an upper bound:
$$
\# \{ p \ \text{prime}:\ y<p\leq x\}\ \leq\ \#\{ n\leq x:\ p\nmid n \ \text{for all primes} \  p\leq y\} .
$$
At the start this works quite well. If $y=2$ the quantity on the right is $\frac 12 x \pm 1$, and so bounded above by $\frac 12 x+1$. If $y=3$ then we remove roughly a third of the remaining integers (leaving two-thirds of them) and so the bound improves to
$\frac 23 \cdot \frac 12 x+2$. For $y=5$ we have four-fifths of the remaining integers to get the upper bound
$\frac 45 \cdot \frac 23 \cdot \frac 12 x+4$. And, in general, we obtain an upper bound of no more than
$$
  \prod_{p\leq y} \left( 1 -\frac 1p \right) \cdot x + 2^{\pi(y)-1}.
$$
The problem with this bound is the second term \dots as one sieves by each consecutive prime, the second term, which comes from a bound on the rounding error, doubles each time, and so quickly becomes larger than $x$ (and thus this is a useless upper bound). This formula   does allow us, by letting $y\to \infty$ slowly with $x$, to prove that
$$
\lim_{x\to \infty} \ \frac{\pi(x)}{x} \ = \ 0;
$$
that is the primes are a vanishing proportion of the integers up to $x$, as $x$ gets larger.\footnote{To deduce this we need to know that $\lim_{y\to \infty} \prod_{p\leq y} \left( 1 -\frac 1p \right) = 0$, a fact established by Euler.}

There has been a lot of deep and difficult work on improving our understanding of the sieve of Eratosthenes, but we are still unable to get a very good upper bound for the number of primes in this way. Moreover we are unable to use the sieve of Eratosthenes (or any other sieve method) to get good {\sl lower bounds} on the number of primes up to $x$

The first big leap in our ability to give good upper and lower bounds on $\pi(x)$ came from an extraordinary observation of Chebyshev in 1851. The observation (as reformulated by Erd\H os in 1933) is that the binomial coefficient $\binom{2n}n$ is an integer, by definition, and is divisible  once  by each prime $p$ in $(n,2n]$, since $p$ is a term in the expansion of the numerator $(2n)!$, but not of the denominator, $n!^2$. Therefore
$$
\prod\Sb p \ \text{prime}\\ n<p\leq 2n\endSb p  \ \leq \ \binom{2n}{n} .
$$
Now, by the binomial theorem,
$\binom{2n}{n} \leq \sum_{j=0}^{2n} \binom{2n}{j}=(1+1)^{2n} = 2^{2n}$.
Moreover for each $p\in (n,2n]$ we have $p>n$ and so
$$
n^{\pi(2n)-\pi(n)} \ = \ \prod\Sb p \ \text{prime}\\ n<p\leq 2n\endSb n \ \leq \  \prod_{n<p\leq 2n} p  \ \leq \ \binom{2n}{n} \leq 2^{2n}.
$$
Taking logarithms this gives us the upper bound $\pi(2n)-\pi(n)\leq \frac{2n\log 2}{\log n}$, and summing this bound yields
$$
\pi(x)\ \leq\ \sum_{j\geq 0} \big(\pi\big(\frac{x}{2^{j}}\big) - \pi\big(\frac{x}{2^{j+1}}\big)\big)\ \leq\ (\log 4+\epsilon) \frac{x}{\log x} ,
$$
for $x$ sufficiently large.

One can interpret this proof as a mixture of algebra and analysis: \ The algebra comes when we consider the primes that divide the central binomial coefficient to get a lower bound on its size; the analysis when we bound the size of the central binomial coefficient by comparing it to the size of other binomial coefficients. This sets the pattern for what is to come.

One can also obtain good lower bounds by studying the prime {\sl power} divisibility of  $\binom{2n}n$:\ First we note that if $p$ divides $\binom{2n}n$ then $p\leq 2n$, since the numerator is the product of all integers $\leq 2n$, and hence can have no prime factor larger than $2n$. The key observation, essentially due to Kummer, is that if a prime power $p^{e_p}$ divides $\binom{2n}n$ then we even have that $p^{e_p}\leq 2n$.
We couple that with the observation that $\binom{2n}n$ is the largest of the binomial coefficients $\binom{2n}j$, and therefore
$$
2^{2n} = (1+1)^{2n} = \sum_{j=0}^{2n} \binom{2n}{j} \leq (2n+1) \binom{2n}n.
$$
Combining this information, we obtain
$$
\frac{2^{2n}}{2n+1}\ \leq\ \binom{2n}n\ =\ \prod\Sb p \ \text{prime}\\
p\leq 2n\endSb p^{e_p}\ \leq \ \prod\Sb p \ \text{prime}\\
p\leq 2n\endSb 2n \ = \  (2n)^{\pi(2n)} .
$$
Taking logarithms gives us the lower bound $\pi(2n) \geq \frac{2n\log 2 - \log(2n+1)}{\log (2n)}$, and therefore
$$
\pi(x) \  \geq \ (\log 2-\epsilon) \frac{x}{\log x}
$$
if $x$ is sufficiently large. Hence we have shown that there exist constants $c_2>1>c_1>0$,\footnote{In fact our proof shows that we can take any $c_1<\log 2$ and   $c_2>\log 4$.} such that if $x$ is sufficiently large then
$$
c_1 \frac{x}{\log x}\ \leq\ \pi(x)\ \leq \ c_2 \frac{x}{\log x}.  \tag{2}
$$
The {\sl Prime Number Theorem}, that is the conjecture of Gauss and Legendre estimating the number of primes up to $x$, can be re-phrased as the claim that these inequalities hold for {\sl any} constants $c_1$ and $c_2$ satisfying $c_2>1>c_1>0$; and in particular we can take both $c_1$ and $c_2$ arbitrarily close to $1$.

Can the method of Chebyshev and Erd\H os be suitably modified to prove the result?  In other words, perhaps we can find some other product of factorials that yields an integer, and in which we can track the divisibility of the large prime factors, so that we obtain constants $c_1$ and $c_2$ in (2), that are closer to $1$. We might expect that the closer the constants get to $1$, the more complicated the product of factorials, and that has been the case in the efforts that researchers have made to date.\footnote{Which is not to say that someone new, not overly influenced by previous, failed attempts, might not come up with a cleverer way to modify the previous approaches.} There is one remarkable identity that gives us hope. First note that in our argument above we might replace $\binom{2n}n$ by $[x]!/[x/2]!^2$, where $x=2n$.\footnote{Here $[t]$ denotes the largest integer $\leq t$.} Hence the correct factorials to consider take the shape $[x/n]!$ for various integers $n$. Our remarkable identity is:
$$
\prod\Sb p \ \text{prime}\\ e\geq 1 \\  p^e\leq x\endSb p \ = \
\prod_{n\leq x}   \left[ x/n \right]  !\ ^{\mu(n)} . \tag{3}
$$
This needs some explanation. The left-hand side is the product over the primes $p\leq x$, each $p$ repeated $k_p$ times, where $p^{k_p}$ is the largest power of $p$ that is $\leq x$. On the right hand side we have the promised factorials, each to the power $-1, 0$ or $1$. Indeed the {\sl M\H obius function} $\mu(n)$ is defined as
$$
\mu(n) \ = \ \cases
0 & \text{if there exists a prime} \ p \ \text{for which} \ p^2 \ \text{divides} \ n;\\
(-1)^k & \text{if} \ n\ \text{is squarefree, and has exactly} \ k \ \text{prime factors}
\endcases
$$
The difficulty with using the identity (3) to prove the prime number theorem is that the length of the product on the right side grows with $x$, so there are too many terms to keep track of.  One idea is to simply take a finite truncation of the right-hand side; that is
$$
\prod_{n\leq N}   \left[ x/n \right]  !\ ^{\mu(n)}
$$
for some fixed $N$. The advantage of this is that, once $x>N^2$, then this product is divisible by every prime $p\in [x/(N+1),x]$ to the power $1$. The disadvantage is that the product is often not an integer, though we can correct that by multiplying through by a few smaller factorials.  We can handle this and   other difficulties that arise, to obtain (2) with other values of $c_1$ and $c_2$, which each appear to be getting closer and closer to $1$. However when we analyze what it will take to prove that the constants (which we now denote by $c_1(N)$ and $c_2(N)$ since they depend on $N$) tend to $1$ as $N\to \infty$, we find that the issue lies in the average of the exponents $\mu(n)$. In fact one can prove that the constants $c_1(N)$ and $c_2(N)$ do tend to $1$ if
$$
\lim_{N\to \infty} \ \frac 1N \ \sum_{n\leq N} \mu(n) \ \ \text{exists and equals} \ 0. \tag{5}
$$
That is, (5)  is  equivalent to the prime number theorem.

The problem in (5) certainly looks more approachable  than the prime number theorem itself, even though the problems are equivalent. It can be restated as: There are roughly the same number of squarefree integers with an even number of prime factors, as there are squarefree integers with an odd number of prime factors. This seems very plausible, and leads to many elementary approaches, as we will discuss in section 9, and beyond.

One of the most famous old problems about primes was  Bertrand's postulate, to prove that for every $n>1$ there is always a prime $p$ for which $n<p<2n$. This follows easily from suitable modifications of the above discussion with the binomial coefficient $\binom{2n}n$.
This beautiful proof, due to Erd\H os at age 20, announced his arrival onto the mathematical stage and inspired the lines:
\smallskip

{\sl Chebyshev said it, and I say it again:

There is always a prime between $n$ and $2n$.}
\medskip

Up to now we have proved that $\pi(x)$ lies between two multiples of $x/\log x$, and we have looked to see whether the ratio of $\pi(x)$ to $x/\log x$ tends to $1$, as predicted by Gauss and Legendre.  Is any other behaviour possible, given what we know already?  There are two possibilities:\ As $x$ grows larger,

\roster
\item"(i)" The ratio of $\pi(x)$ to $x/\log x$ oscillates, never tending to a limit.
\item"(ii)" The ratio of $\pi(x)$ to $x/\log x$ tends to a limit, but that limit is not $1$.
\endroster
\medskip

Our goal in the rest of this section is to show that option (ii) is not possible. Indeed we will show that if there is a limit then that limit would have to be $1$. Yet again the trick is to study factorials  both algebraically (by determining their prime factors), and analytically (by analyzing their size).

By definition,
$$
\log N! \ = \ \sum_{n=1}^N \log n.
$$
The right-hand side is very close to the integral of $\log t$ over the same range. To see this note that the logarithm function is monotone increasing, which  implies that
$$
 \int_{n-1}^n \log t\ dt < \log n < \int_{n}^{n+1} \log t\ dt
$$
for every $n\geq 1$. Summing these inequalities  over all integers $n$ in the range $2\leq n\leq N$ (since $\log 1=0$), we obtain that for $N\geq 2$, the value of $\log N! $ equals
$$
\int_1^N \log t\ dt \ = \ [t(\log t-1)]_1^n \ = \ N(\log N-1)+1 ,
$$
plus an error that is no larger, in absolute value, than $\log N$.

On the other hand $N!$ is the product of the integers up to $N$, and we want to know how often each prime divides this product.  The integers $\leq N$ that are multiples of a given integer $m$ (which could be a prime or prime power) are $m, 2m,\ldots, [N/m]m$, since $[N/m]m$ is the largest multiple of $m$ that is $\leq N$, and therefore there are $[N/m]$ such multiples.  Now the power of prime $p$ dividing $N!$ is given by the number of integers $\leq N$ that are multiples of $p$, plus the number of integers $\leq N$ that are multiples of $p^2$, etc., which yields a total of
$$
 \left[ \frac{N}{p} \right] + \left[ \frac{N}{p^2} \right] + \left[ \frac{N}{p^3} \right] +\ldots 
$$
Hence, by studying the prime power divisors of $N!$ we deduce that
$$
\log N! \ = \ \sum\Sb p \ \text{prime}\\  p \leq N\endSb \ \log p \left(
\left[ \frac{N}{p} \right] + \left[ \frac{N}{p^2} \right] + \ldots +\right) .
$$
The total error created by discarding those  $[N/p^k]$ terms with $k\geq 2$, and by replacing each $[N/p]$ by $N/p$, adds up to no more than a constant times $N$. Comparing our two estimates for $\log N!$, and dividing through by $N$ we deduce that
there exists a constant $C$ for which
$$
\left| \sum_{p\leq N} \frac{\log p} p \ -\ \log N \right| \ \leq \ C \tag{7}
$$
for all $N\geq 1$.

Now let's suppose that there exists a constant $\eta$, such that for any $\epsilon>0$
$$
(\eta-\epsilon) \frac{x}{\log x}\ \leq\ \pi(x)\ \leq \ (\eta+\epsilon) \frac{x}{\log x},
$$
for all sufficiently large $x$ (say $>x_\epsilon$). Our goal is to show that $\eta=1$.

We will work with the following identity:
$$\align
\sum_{p\leq N} \frac{\log p} p \ &=  \
\sum_{p\leq N}  \left\{   \frac{\log N}{N} + \int_p^N \frac{\log x-1}{x^2} \ dx  \right\}  \\
&=\ \pi(N) \frac{\log N}{N} + \int_2^N \pi(x)\frac{\log x-1}{x^2} \ dx  ,
\endalign
$$
inserting our assumed bounds on $\pi(x)$ to obtain bounds on $\sum_{p\leq N} \frac{\log p} p$. The part of the integral with $x\leq x_\epsilon$ is bounded by some constant that only depends on $\epsilon$, call it $C_1(\epsilon)$. Therefore we obtain an upper bound
$$\align
\sum_{p\leq N} \frac{\log p} p \ &\leq C_1(\epsilon) +  \pi(N) \frac{\log N}{N} + \int_{x_\epsilon}^N \pi(x)\frac{\log x-1}{x^2} \ dx \\
&\leq C_1(\epsilon) +  (\eta+\epsilon) \frac{N}{\log N} \frac{\log N}{N} + \int_{x_\epsilon}^N (\eta+\epsilon) \frac{x}{\log x} \frac{\log x-1}{x^2} \ dx \\
&\leq C_2(\epsilon) + (\eta+\epsilon)  \int_{x_\epsilon}^N \frac{dx}{x} \ \leq\  (\eta+\epsilon) \log N + C_2(\epsilon),
\endalign
$$
for some constant $C_2(\epsilon)$, that only depends on $\epsilon$.  This implies that $\eta\geq 1$ else we let $\epsilon=(1-\eta)/2$ and this inequality contradicts (7) for $N$ sufficiently large. An analogous proof with the lower bound implies that $\eta\leq 1$. We deduce that   if $\pi(x) \big/ \frac x{\log x}$ tends to a limit as $x\to \infty$ then that limit must be $1$.

\subhead 3. A first reformulation:\ Introducing appropriate weights \endsubhead So far, we have counted primes by estimating the size of the product of the primes in some interval. Taking logs, this means that we bounded
$$
\sum\Sb p \ \text{prime}\\  p\leq x\endSb  \log p.
$$
We denote this by $\theta(x)$; and we also define its close cousin,\footnote{The reader might verify  that $\theta(x)$ and $\psi(x)$ do not differ by more than a bounded multiple of $\sqrt{x}$.}
$$
\psi(x):= \ \sum\Sb p \ \text{prime}\\ m\geq 1 \\ p^m\leq x\endSb  \log p.
$$
We will see that when we do calculations, these functions seem to be more natural than $\pi(x)$ itself. This fits rather well with Gauss's original musings in his letter to Encke. The key phrase is:

\block {\eightpoint
    {\sl I soon recognized, that under all variations of this frequency {\rm [of prime numbers]}, on  average, it is nearly inversely proportional to the logarithm.}  }
    \endblock
\smallskip

\noindent We re-word this as ``{\sl The density of primes at around $x$ is $1/\log x$.}'' Then we would expect that the number of primes, each weighted with $\log p$  (that is, the sum $\theta(x)$) should be well-approximated by
$$
\int_2^x \log t  \cdot \frac{dt}{\log t} \ = \ \int_2^x dt = x-2.
$$
Occam's razor tells us that, given two choices, one should opt for the more elegant one. There can be little question that
$x$ is a   more pleasant function to work with than the complicated integral $\Li(x)$, and so we will develop the theory with logarithmic weights, and therefore use the function $\theta(x)$ rather than $\pi(x)$.\footnote{Using partial summation, it is not difficult to show that a good estimate for one is equivalent to an analogous estimate for the other, so there is no harm done in focusing on $\theta(x)$.}

We believe that
$$
\left|  \sum_{p\leq x} \log p - x \right|  \leq  x^{1/2} (\log x)^2 ,\tag{RH2}
$$
since this is equivalent to our conjecture (RH1) on $\pi(x)-\Li(x)$. The prime number theorem is equivalent to the much weaker assertion that
$$
\lim_{x\to \infty} \ \frac 1x \sum_{p\leq x} \log p \ \ \text{exists and equals} \ 1.
$$

\subhead 4. Riemann's memoir \endsubhead
In a nine page memoir written in 1859, Riemann outlined an extraordinary plan to  attack the  elementary question of counting prime numbers using deep ideas from the theory of complex functions. His approach begins with what we now call the {\sl Riemann zeta-function}:
$$\zeta(s) :=\sum_{n\geq 1} \frac 1{n^s}.$$
To make sense of an infinite sum it needs to converge, and preferably be absolutely convergent, meaning that we can re-arrange the order of the terms without changing the value.\footnote{Riemann proved that if one has a convergent but not absolutely convergent sum then one might get different limits if one re-arranges the order of the terms in the sum.} The sum defining $\zeta(s)$ is absolutely convergent only when Re$(s)>1$.  This is especially interesting when we apply the {\sl Fundamental Theorem of Arithmetic} to each term in the sum: Every integer $n\geq 1$ can be factored in a unique way  and, vice-versa, every product of primes yields a unique positive integer. Then we can write
$$
n=2^{n_2}3^{n_3}\ldots ,
$$
where each $n_j$ is a non-negative integer, and only finitely many of them are non-zero. Hence
$$\align
\zeta(s) \ & = \ \sum\Sb n_2, n_3, n_5,\ldots \geq 0\endSb \frac 1{(2^{n_2}3^{n_3}5^{n_5}\ldots)^s} \\
&= \prod\Sb p\ \text{prime} \endSb \left( \sum_{n_p\geq 0} \frac 1{(p^{n_p})^s} \right)  \ = \
 \prod\Sb p\ \text{prime} \endSb \left( 1 - \frac 1{p^s}\right)^{-1} .
 \endalign
$$
This product over primes is an {\sl Euler product} --- indeed, it was Euler who first seriously explored the connection between $\zeta(s)$ and the distribution of prime numbers, though he did not penetrate the subject as deeply as Riemann.

The Euler product provides a connection between $\zeta(s)$ and prime numbers, and this was exploited by Riemann in an interesting way. Since the sum defining $\zeta(s)$ is absolutely convergent when Re$(s)>1$, it is safe to perform calculus operations on $\zeta(s)$ in this domain. By taking the logarithmic derivative we have
$$
-\frac{\zeta'(s)}{\zeta(s)} \ = \ -\frac{d}{ds} \log \zeta(s)
\ = \  \sum\Sb p\ \text{prime} \endSb \frac{d}{ds}  \log \left( 1 - \frac 1{p^s}\right)
$$
using the Euler product. Using the chain rule, we then obtain
$$
-\frac{\zeta'(s)}{\zeta(s)} \ =  \  \sum\Sb p\ \text{prime} \endSb  \frac{  \log p } {p^s-1} \ =  \  \sum\Sb p \ \text{prime}\\ m\geq 1\endSb  \frac{  \log p } {p^{ms}} ;
$$
and notice that the sum of the coefficients of $1/n^s$ on the right-side, for $n$ up to $x$, equals $\sum_{p^m\leq x} \log p = \psi(x)$.  As we remarked above, this is a close cousin of $\theta(x)$, and we have now observed, like Riemann, that it arises naturally in this context.

\subhead 5. Contour integration \endsubhead One of the great discoveries of 19th century mathematics is that it is possible to convert problems of a discrete flavour, in number theory and combinatorics, into questions of complex analysis. The key lies in finding suitable analytic identities to describe combinatorial issues. For example, if we ask whether two integers, $a$ and $b$, are equal, that is ``does $a=b$?'', then this is equivalent to asking ``is $a-b=0$?'' and therefore we need some analytic device that will distinguish $0$ from all other integers. This is given by the integral of the exponential function around the circle:
$$
\int_0^1 e^{2i\pi nt} dt = \cases 1 & \text{if} \ n=0, \\
0 & \text{otherwise} .
\endcases
$$
So, for example, if we want to determine the number of pairs $p,q$ of primes  $\leq n$, which add to give the even integer $n$, we create an integral that gives $1$ if $p+q=n$, and $0$ otherwise, and then sum over all such $p$ and $q$. Therefore we have
$$
\#\{ p,q\leq n:\ p,q \ \text{primes,} \ p+q=n\} \ = \
\sum\Sb p,q\leq n\\ p,q \ \text{primes} \endSb    \int_0^1 e^{2i\pi (p+q-n)t} dt
$$
and this can be re-arranged as
$$
\int_0^1 e^{-2i\pi nt}  \left(  \sum\Sb p \leq n\\ p \ \text{prime } \endSb    \  e^{2i\pi pt} \right)^2            dt .
$$
This is of course an approach to {\sl Goldbach's conjecture} (that every even integer $\geq 4$ is the sum of two primes), and it is (arguably) surprising that understanding this integral is equivalent to the original combinatorial number theory question.

So we have seen how to analytically identify when two integers are equal, and why that is useful. Next we will show how to analytically verify a proposed inequality between two real numbers. As you might have guessed, we start by noting that asking whether $u<v$ is the same as asking whether $v-u>0$, and so we restrict our attention to determining whether a given real number is $<0$ or $=0$ or $>0$ (though we are less interested in the middle case). Here the trick is that for any $\sigma>0$ we have
$$
\frac 1{2 \pi} \int_{-\infty}^\infty \frac{e^{y(\sigma+it)}}{\sigma+it}\ dt \ = \
\cases 0 &\text{if} \ y<0; \\  1/2  &\text{if} \ y=0; \\  1 &\text{if} \ y>0,\endcases
$$
which is  {\sl Perron's formula}. It is convenient to write $s$ for $\sigma+it$ and, instead of the limits of the integral, we write  ``Re$(s)=\sigma$'', understanding that we take $s$ along the line Re$(s)=\sigma$, that is,
$s=\sigma+it$ as $t$ runs from $-\infty$ to $\infty$. Hence we integrate $e^{ys}/s$.
Moreover if we let $z=e^y$, the formula can be rephrased as
$$
\frac 1{2i\pi} \int_{\text{Re}(s)=\sigma}  \frac{z^{s}}{s}\ ds \ = \
\cases 0 &\text{if} \ 0<z<1; \\  1/2  &\text{if} \ z=1; \\  1 &\text{if} \ z>1,\endcases
$$
for any $\sigma>0$.
In number theory, the most common use of Perron's formula is to identify when an integer $n$ is $< x$, that is
when $x/n>1$.

We are interested in estimating $\psi(x)=\sum_{p^m\leq x} \log p$. We extend the sum to all prime powers, multiplying by $1$ if $p^m\leq x$, and by $0$ otherwise, which we achieve by using Perron's formula with $z=x/p^m$. The outcome is $\psi^*(x)$ which has the same value as $\psi(x)$ except that we subtract $\frac 12 \log x$ if $x$ is a prime power. Therefore we have
$$
\psi^*(x)\ = \ \sum\Sb p \ \text{prime} \\ m\geq 1\endSb \log p \cdot
\frac 1{2i\pi} \int_{\text{Re}(s)=\sigma}  \frac{(x/p^m)^{s}}{s}\ ds ,
$$
for any $\sigma>0$. Now we would like to swap the order of the summation and the integral, but there are convergence issues. Fortunately these are easily dealt with when the sum is absolutely convergent, as happens when $\sigma>1$. Then we have, after a little re-arrangement,
$$\align
\psi^*(x)\ &= \ \frac 1{2i\pi} \int_{\text{Re}(s)=\sigma}   \sum\Sb p \ \text{prime} \\ m\geq 1\endSb \frac{\log p}{p^{ms}} \cdot
 \frac{x^{s}}{s}\ ds \\
 &=\ \frac 1{2i\pi} \int_{\text{Re}(s)=\sigma}    - \frac{\zeta'(s)}{\zeta(s)} \cdot  \frac{x^{s}}{s}\ ds . \tag{11}
 \endalign
$$
This seems, at first sight, to be a rather strange thing to do. We have gone from a perfectly understandable question like estimating $\pi(x)$, involving a sum that is easily interpreted, to a rather complicated integral, over an infinitely long line in the complex plane, of a function that is delicate to work with in that it is only well-defined when Re$(s)>1$. It is by no means obvious how to proceed from here, as we will discuss in more detail in the next section.

The proof of (11) did not use many properties of $\zeta'(s)/\zeta(s)$. In fact if $a_n$ is any sequence of real numbers with each $|a_n|\leq 1$ then define the {\sl Dirichlet series} $A(s)=\sum_{n\geq 1} a_n/n^s$, to obtain
$$
\sum_{n\leq x} a_n \ = \  \frac 1{2i\pi} \int_{\text{Re}(s)=\sigma}      A(s)  \cdot  \frac{x^{s}}{s}\ ds .
$$

\subhead{6. Riemann's genius} \endsubhead How do  we evaluate the integral in (11)? In complex integration the idea is to shift the path of the integral to one on which the integrand is ``very small''. Then the value of the integral is given by the sum of the ``residues'' of the integrand at its ``poles.'' There is a lot to explain here --- indeed the main points of a first course in complex analysis. Rather than get into all of these details, let me just say that the poles are the points where the function goes to $\infty$, like the point $s=1$ for the function $1/(s-1)$. And if the function $f(s)$ has a pole at, say, $s=1$, whereas $(s-1)f(s)$ equals $r\in \Bbb C$ at $s=1$, then we say that $f(s)$ has a {\sl simple pole  at $s=1$ with residue $r$}.

So what new path should we take from $\sigma-i\infty$ to $\sigma+i\infty$ to be able to apply this strategy to the integral in (11)? If we are going to choose the same path for each value of $x$, then we want to ensure that the integrand (on the path) does not grow large with $x$. Now $|x^s|=x^{\text{Re(s)}}=x^\sigma$, so the smaller $\sigma$ is, the better. In fact if we make $\sigma$ negative then the $x^s$ in the integrand will ensure that the integral over this line gets smaller as $x\to \infty$. Or, rather, that would be true if  $\zeta(s)$  and $(\zeta'/\zeta)(s)$ are defined in this region, which they are not, for now.

Under the right conditions, functions that are naturally defined only in part of the complex plane (like $\zeta(s)$), can be re-defined so that the function can be appropriately extended to the rest of the complex plane. This is called an {\sl analytic continuation}, which involves what can be a deep and subtle theory. In such circumstances, one can express the function in terms of a Taylor series,\footnote{That is, for $f(s)$ at $s_0\in \Bbb C$, there exist constants $a_0,a_1,\ldots$ and some  constant $r$, such that if $|s-s_0|<r$ then
$a_0+a_1(s-s_0)+a_2(s-s_0)^2+\ldots$ is absolutely convergent, and converges to $f(s)$.}
or a Laurent series if there is a pole.
Analytic continuations are a little bit mysterious -- for example the theory allows for more than one apparently different way that one can define the function on the rest of the complex plane, but it will turn out that any two such definitions will have equal values everywhere they are both defined.\footnote{This allows us to give the analytic continuation the same name as the original function, since we know that it can only be analytically continued in one way, if at all.} Anyway, we can analytically continue $\zeta(s)$ to all of the complex plane, except for its pole at $s=1$, and to do this, Riemann discovered some remarkable properties of $\zeta(s)$ (which we do not pursue here). There are several subtleties involved in bounding the contribution of the integrand on the new contours, and Riemann succeeded in doing that. Finally one needs to find  the poles of
$$ - \frac{\zeta'(s)}{\zeta(s)} \cdot  \frac{x^{s}}{s},$$
and to compute their residues:\
Evidently $x^s$ has no poles in the complex plane, and $1/s$ has a simple pole at $s=0$, which contributes a residue of
$$
- \lim_{s\to 0}\  s \ \frac{\zeta'(s)}{\zeta(s)} \  \frac{x^s}s =
 \lim_{s\to 0}\    \frac{\zeta'(s)}{\zeta(s)} \  x^s \ = \ - \frac{\zeta'(0)}{\zeta(0)} \   x^{0} \ = \ - \frac{\zeta'(0)}{\zeta(0)}
$$
to the value of the integral.
The poles of $\zeta'(s)/\zeta(s)$ are the poles and zeros of $\zeta(s)$. The only pole of $\zeta(s)$ is at $s=1$, and so this contributes the  residue
$$
-\lim_{s\to 1}\  (s-1)\frac{\zeta'(s)}{\zeta(s)} \  \frac{x^s}s  \  = \
-\lim_{s\to 1}\  (s-1)\left(\frac{-1}{(s-1)}\right) \  \frac{x^1}1  \  = \ x,
$$
the expected main term, to the value of the integral (since $\zeta(0)\ne 0$).   The Euler product representation of $\zeta(s)$ converges in Re$(s)>1$,  so there can be no zeros of $\zeta(s)$ in this half-plane. Otherwise the zeros of $\zeta(s)$ are rather mysterious. All we can really say is that if $\zeta(s)$ looks like $c(s-\rho)^m$, near to the zero $\rho$, for some integer $m\geq 1$ (and non-zero constant $c$) then
$-\zeta'(s)/\zeta(s)$ looks like $-m/(s-\rho)$, and therefore the residue at $s=\rho$
is
$$
-\lim_{s\to \rho}\  (s-\rho)\frac{\zeta'(s)}{\zeta(s)} \  \frac{x^s}s \  = \
-\lim_{s\to \rho}\  (s-\rho)\left(\frac{m}{(s-\rho)}\right) \  \frac{x^\rho}\rho   \  = \ -m \ \frac{x^\rho}\rho.
$$
If we count such a zero of multiplicity $m$, $m$ times in the sum, then we have evaluated the integral so as to yield
Riemann's remarkable {\sl explicit formula}:
$$
\boxed{
\qquad \psi^*(x)\ =\ x  - \sum_{\rho:\ \zeta(\rho)=0} \frac{x^\rho}\rho -\frac{\zeta'(0)}{\zeta(0)} . \qquad }
$$
The left-hand side is a step function, with a jump at each prime power, whereas the right-hand side is the sum of infinitely many smooth functions. Somehow these smooth functions, which correspond to the zeros of $\zeta(s)$ conspire to stay constant as $x$ varies, other than to jump at exactly the prime powers.
Surprising, it may be, but is it useful? We have gone from a simple question like counting the number of primes up to $x$, to a sum over all of the zeros of the analytic continuation of $\zeta(s)$.  Ever since Riemann's memoir, mathematical researchers have struggled to find a way to fully understand the zeros of $\zeta(s)$, so as to make this ``explicit formula'' useful. We have had some, rather limited, success.

Riemann showed that there are infinitely many zeros of $\zeta(s)$ so we have a problem in that the sum over $\rho$, in Riemann's explicit formula, is an infinite sum and one can easily show that it is not absolutely convergent. So to evaluate it directly, we would need to detect cancellation amongst the summands, something that we are not very skilled at. Instead, one can modify Riemann's argument to show that one can truncate the sum, taking only those $\rho$ in the box  up to height $T$,
$$
B(T):=\{ \rho\in \Bbb C:\ 0\leq \text{Re}(\rho) \leq 1 , \ \ -T\leq \text{Im}(\rho) \leq T\} ,
$$
in our sum. This turns out to be a finite set, and we get the  explicit formula,
$$
\psi(x)\ =\ x  - \sum\Sb \rho:\ \zeta(\rho)=0 \\  \rho\in B(T)\endSb \frac{x^\rho}\rho \ +  \ \text{a small error},
$$
where the ``small error'' is small if $T$ is appropriately chosen (as a function of $x$; typically $T=\sqrt{x}$). Then we can  bound $|\psi(x)-x|$, by taking absolute values in the sum over zeros $\rho$, above.  Again the key issue is that $|x^\rho|=x^{\text{Re}(\rho)}$ so that
$$
\left| \sum\Sb\rho:\ \zeta(\rho)=0 \\ \rho\in B(T)\endSb \frac{x^\rho}\rho \right| \ \leq \ \sum\Sb\rho:\ \zeta(\rho)=0 \\ \rho\in B(T)\endSb \left| \frac{x^\rho}\rho \right| \ =\ \sum\Sb\rho:\ \zeta(\rho)=0 \\ \rho\in B(T)\endSb  \frac{x^{\text{Re}(\rho)}}{|\rho|}   \ \leq \ x^{m(T)} \sum\Sb\rho:\ \zeta(\rho)=0 \\ \rho\in B(T)\endSb \frac{1}{|\rho|}
$$
where 
$$m(T):=\max_{\rho\in B(T)} \text{Re}(\rho).$$
The sum over zeros can be shown to be bounded by a multiple of $(\log T)^2$, so if we can get a good bound on $m(T)$ then we will be able to deduce the prime number theorem. By ``good bound'' here we mean that $m(T)$ must be somewhat less than $1$, in fact
$$ m(T) \leq 1 - \frac{3\log\log T}{\log T} $$
will do.

Riemann  made a few calculations of the zeros of $\zeta(s)$ and all the  real parts  seemed to be $1/2$. This led to him to:\footnote{Riemann actually wrote: ``{\sl It is very probable that all roots are {\rm [on the $\frac 12$-line]}. Certainly one would wish for a stricter proof here; I have meanwhile temporarily put aside the search for this after some fleeting futile attempts, as it appears unnecessary for the next objective of my investigation.}''}

\proclaim{The Riemann Hypothesis} If $\zeta(\rho)=0$ with
$0\leq \text{\rm Re}(\rho)\leq 1$, then $ \text{\rm Re}(\rho)=\frac 12$.
\endproclaim

Even though the Riemann Hypothesis remains unproven today, more than 150 years after Riemann's article, most mathematicians believe that it is true. There have been extensive calculations, proving that the first ten billion  zeros above the real axes  lie on the $\frac 12$-line.  More persuasive is that it fits so well with so many other ideas that the world would be a much uglier place if it is not true. The reason that we are drawn to the $\frac 12$-line is Riemann's remarkable {\sl functional equation} which shows that $\zeta(s)$ can easily be determined in terms of $\zeta(1-s)$; in particular once we understand $\zeta(s)$ for Re$(s)\geq \frac 12$ then we understand it on the whole complex plane.

If, the Riemann Hypothesis is true  then each $x^\rho$ has absolute value
$x^{1/2}$ or $<1$, and one can deduce, via the argument that we have just sketched, the estimates (RH1) and (RH2) for $\pi(x)$ and $\theta(x)$. Actually there is a very intimate link between upper bounds for the real parts of the zeros of $\zeta(s)$ and bounds on the error term in the prime number theorem, and one can show that if either (RH1) or (RH2) is true then the Riemann Hypothesis follows. This connection goes much further.  For example, fix $1>\beta>1/2$. Then  all zeros of $\zeta(s)$ satisfy
\smallskip
\centerline{ Re$(s)<\beta$ if and only if  $\left| \sum_{p \leq x} \log p - x \right| \leq C_\beta x^{\beta}$,}
\smallskip

\noindent for some constant $C_\beta>0$.  How strange! Here we are in two different worlds, counting primes, and zeros of the analytic continuation of a function, and yet a key part of understanding each is equivalent. This is the bedrock on which mathematics is formed. Surprising connections between fields that have no obvious right to be related, and yet they are, at some fundamental level.  Riemann's work gave one of the first results of this type, and now every field of research in pure mathematics is full of such links.

Riemann's connection is not restricted to this one question. Indeed, using the explicit formula, one can {\sl reformulate} many different problems about primes, as problems about zeros of  zeta-functions, upon which we can use the tools of analysis.  Mathematicians love  bringing  fields together that seem  so distant,   hopefully allowing a more rounded perspective of both.

These observations are so seductive that they have been the thrust of almost all research into the distribution of prime numbers ever since.  Moreover there are  many other good questions about prime numbers, number fields, finite fields, curves and varieties that  can be re-cast in terms of appropriate zeta-functions, so there is no end to what can be investigated by such methods.

\subhead 7. The coup de gr\^ace in the proof of the prime number theorem\endsubhead
What we have sketched above is not quite the end of the story of the proof of the prime number theorem. Although Riemann came up with the
whole idea, and made many spectacular advances in his short memoir, he could not
give an unconditional proof. He left several steps to be completed. These turned out to be very difficult indeed, and  it was only 37 years later that Hadamard and de le Vall\'ee Poussin did so, independently, in 1896. The final step that was left to them, was to show that $\zeta(s)$ has no zeros on the line Re$(s)=1$, which we call the {\sl $1$-line} from now on. Both of their proofs, and most that followed, show that $\zeta(s)$ cannot have a zero at $1+it$ by showing: 
$$
\text{If} \ \ \zeta(1+it)=0 \ \ \text{then} \ \ \zeta(1+2it)=\infty; \tag{13}
$$ 
that is, $\zeta(s)$ has a pole at $1+2it$. However we have already noted that $\zeta(s)$ only has a pole at $s=1$, and therefore $t=0$. This yields  a contradiction to the assumption that $\zeta(1+it)=0$. The proofs of Hadamard and de le Vall\'ee Poussin are complicated, and the proof of Mertens that can be found in every textbook is relatively easy without being enlightening. Nonetheless it is not difficult to get an intuitive  feel for why (13) should be true: \ Since $\zeta(s)$ is an analytic function, if it  equals $0$ at $1+it$, then it must be well-approximated by the leading term in its Taylor series, $c(s-(1+it))^r$, when $s$ is sufficiently close to $1+it$, for some integer $r\geq 1$ and some non-zero constant $c$. For example,

\centerline{ if  $s=1+it+\frac 1{\log x}$ then
$\zeta(s)\approx c/(\log x)^r$,}

\noindent which is pretty small.\footnote{Throughout we will use the notation ``$\approx$'' to mean ``is approximately equal to''. Typically I'll avoid being too precise as this can introduce uninteresting yet substantial technicalities.}

Since Re$(s)>1$ we can determine $\zeta(s)$ in terms of its Euler product, and one can approximate $\zeta(s)$ well  at $s=1+it+\frac 1{\log x}$  by truncating the Euler product at $x$. In other words
$$
\zeta(s) \ \approx \ \prod\Sb p\ \text{prime} \\ p\leq x \endSb \left( 1 - \frac 1{p^{1+it}}\right)^{-1} .
$$
Now the $p$th term in this Euler product has absolute value
$$\big|1 - \frac 1{p^{1+it}}\big|^{-1}\geq \big(1+\frac 1p\big)^{-1} , \tag{17} $$
and one can deduce from (7) that $$
\prod\Sb p\ \text{prime} \\ p\leq x \endSb \left( 1 + \frac 1{p}\right)^{-1}\ \approx \ \frac {c'}{\log x},
$$
for some constant $c'$.  This implies the lower bound $|\zeta(s)|\geq c''/\log x$.

Comparing these two estimates for $|\zeta(s)|$ allows us to deduce that
$$
c/(\log x)^r \geq c'''/\log x$$ for all sufficiently large $x$. Hence $r\leq 1$, but we know that $r$ is an integer $\geq 1$, and so $r=1$. Therefore   (17) must be close to equality, most of the time;  that is $-1/p^{1+it} \ \approx \ 1/p$, and therefore
$$p^{it}\approx -1$$
for ``most primes $p$'', a concept we will make precise in section 9.  In fact we will   deduce this directly from (5) by more elementary methods.

Squaring $p^{it}\approx -1$, one sees that $p^{2it}\approx (-1)^2=1$ for most primes $p$, which tells us that if $s=1+2it+\frac 1{\log x}$ then
$$
\zeta(s) \ \approx \  \ \prod\Sb p\ \text{prime} \\ p\leq x \endSb \left( 1 - \frac 1{p^{1+2it}}\right)^{-1}  \ \approx \  \ \prod\Sb p\ \text{prime} \\ p\leq x \endSb \left( 1 - \frac 1{p}\right)^{-1} \approx c'' \log x,
$$
and hence, letting $x\to \infty$, we deduce that  $\zeta(s)$ has a pole at $s=1+2it$.

This is not the only way to show that we cannot have $p^{it}\approx -1$ for most primes  $p$. My favourite technique is to take logarithms and to show that if $p^{it}\approx -1$ for most primes  $p$ then these primes are clustered in intervals of the form
$$
 [ (1-\epsilon)e^{i\pi (2k+1)/|t|} , \ (1+\epsilon)e^{i\pi (2k+1)/|t|} ]
$$
where $k$ is an integer, and $\epsilon$ is very small. One can then use sieve techniques (the direct descendants of the sieve of Eratosthenes, specifically the Brun-Titchmarsh Theorem) to show that primes cannot be clustered into intervals at more than double the expected density, and thus we obtain a contradiction.

Later we will study (5) to   show that ``$p^{it}\approx -1=\mu(p)$ for most primes  $p$'' is impossible.

\subhead 8. Selberg's elementary approach   \endsubhead The explicit formula which directly relates the primes to the zeros of $\zeta(s)$, suggests a tautology between primes and zeros. This persuaded no lesser authorities than Hardy, Ingham and Bohr to assert that it would be impossible to find an {\sl elementary proof} of the prime number theorem.
After all, how could it be possible? The prime number theorem implies restrictions on the zeros of the analytic continuation  of $\zeta(s)$ -- how could one have a proof of that which does not use analysis? As Hardy said in Copenhagen in 1921:

\block {\eightpoint
    {\sl No elementary proof of the prime number theorem is known, and one may
ask whether it is reasonable to expect one. Now we know that the theorem
is roughly equivalent to ... the theorem
that Riemann's zeta function has no roots on a certain line. A proof of such
a theorem, not fundamentally dependent on the theory of functions, seems to
me extraordinarily unlikely. It is rash to assert that a mathematical theorem
cannot be proved in a particular way; but one thing seems quite clear. We have
certain views about the logic of the theory; we think that some theorems ... `lie deep' and others nearer to the surface. If anyone produces an elementary
proof of the prime number theorem, [s]he will show that these views are wrong,
that the subject does not hang together in the way we have supposed, and thatf
it is time for the books to be cast aside and for the theory to be rewritten.}}
\endblock

\noindent And, as Ingham wrote in the introduction of his 1932 book [I1]:

\block {\eightpoint
    {\sl Every known proof of the prime number theorem is based on a certain property of the complex zeros of $\zeta(s)$, and this conversely is a simple consequence of the prime number theorem itself. It seems therefore clear that this property must be used (explicitly or implicitly) in any proof based on $\zeta(s)$, and it is not easy to see how this is to be done if we take account only of real values of $s$.}}
\endblock

The key to Selberg's elementary approach\footnote{I asked Selberg, in around 1989, how he would define ``elementary''. He responded that there is no good definition, but that it is perhaps best expressed as ``what a good high school student could follow.''}   is Selberg's formula:
$$
\log x \sum\Sb  p\leq x  \\  p \ \text{prime} \endSb \log p\ + \
\sum\Sb  pq\leq x  \\  p, \ q \ \text{both prime}  \endSb \log p \
\log q \ =  \ 2x\log x + \text{Error}. \tag{19}
$$
Here the ``Error'' term is bounded by a multiple of $x$. So instead of getting an accurate estimate for the weighted number of primes up to $x$, Selberg gets an accurate estimate for the weighted number of primes and P2's up to $x$ (where a ``P2'' is an integer that is  the product of two primes). Moreover
Selberg [S2] gave an elementary proof that (19) is true using combinatorial methods; and it is tempting to believe that it should not then be difficult to remove the P2's from the equation. But first we ask, how can a formula like (19) hold without any hint of the zeros of $\zeta(s)$?

Selberg does not indicate how he came up with such a formula, and why he would have guessed that it would be true, so we can only speculate. Selberg was a master analyst, so it is plausible that he reasoned as follows:\ The main problem in using Riemann's formula is that if a zero, $\rho$, of $\zeta(s)$ has real part equal to $1$ (which is not easy to disprove), then the corresponding error term, $x^\rho/\rho$, has size $cx$ for some non-zero constant $c$, a positive fraction of the main term.  So, can we come up with a formula, for a quantity similar to the primes, where one such ``bad zero'' cannot have such a damning effect?

One way to approach this is to try to produce an integrand that is similar to the one that Riemann worked with, but for which there is a double pole at $s=1$, and no new higher order poles elsewhere. The simple way to get a double pole is from the function $\zeta''(s)/\zeta(s)$. This also has the feature that if $\rho$ is a simple zero of $\zeta(s)$ then it is a simple pole of $\zeta''(s)/\zeta(s)$. This is not the case with double or higher order zeros, but we expect them to be rare. Hence if we consider the integral
$$
\frac 1{2i\pi} \int_{\text{Re}(s)=\sigma}      \frac{\zeta''(s)}{\zeta(s)} \cdot  \frac{x^{s}}{s}\ ds
$$
then we have the double pole at $s=1$. One can  compute the residue, using the Taylor expansion, at $s=1$, to be
%s of the various terms as a function $s-1$. Now $(s-1)\zeta(s)=1 +\gamma(s-1)+\ldots$ from which one can deduce that
%$(s-1)^2 \zeta''(s)/\zeta(s)= 2 -2\gamma(s-1)+\ldots$. We also have $x^s=x(1+(s-1)\log x+\ldots)$ and $1/s=1-(s-1)+\ldots$, and so
% $$ (s-1)^2  \frac{\zeta''(s)}{\zeta(s)} \cdot  \frac{x^{s}}{s}  = 2x \left( \frac{1}{(s-1)^2}+\frac{\log x-1-\gamma}{s-1}+\ldots \right) , $$
%and therefore the residue at $s=1$ is
$2x(\log x-1-\gamma)$. If $\rho$ is a simple zero of $\zeta(s)$ then its residue is $c_\rho x^\rho$, for some constant $c_\rho$, and with a bit of luck, all of these on the $1$-line will sum up to no more than a constant times $x$.  In other words, one might guess that the above integral should equal $2x\log x$ plus an error which is bounded by at most some multiple of $x$, even if there are zeros on the $1$-line. Evaluating the integral is tricky in its current form, but once we note that
$$
\frac{\zeta''(s)}{\zeta(s)} = \left( \frac{\zeta'(s)}{\zeta(s)} \right)' + \left( \frac{\zeta'(s)}{\zeta(s)} \right)^2,
$$
we can rewrite the integral as
$$
\frac 1{2i\pi} \int_{\text{Re}(s)=\sigma}      \left( \frac{\zeta'(s)}{\zeta(s)} \right)' \frac{x^{s}}{s}\ ds +
\frac 1{2i\pi} \int_{\text{Re}(s)=\sigma}     \left( \frac{\zeta'(s)}{\zeta(s)} \right)^2  \frac{x^{s}}{s}\ ds
$$
which, by Perron's formula, equals
$$
\sum\Sb p^m\leq x \\ p \ \text{prime} \\ m\geq 1\endSb m(\log p)^2 \ + \
\sum\Sb pq\leq x \\ p, q \ \text{primes}  \endSb \log p \log q.
$$
It is easy to show that the prime powers do not contribute much to the first sum, and that $\log p$ is close to $\log x$ for most of the primes $p$ counted in the sum. Hence we get the left-hand side of (19).

This is perhaps why Selberg believed that something like (19) holds, and why it should be accessible to an elementary proof.  However to make this argument elementary required substantial ingenuity (see  [S2]).

How can we deduce the prime number theorem from (19)? The first thing to do is to recast this in terms of the function $\theta(x)$  or, even better, the error term $E(x):=\theta(x)-x$. Using (7) and (19)  we obtain
$$
E(x) \log x + \sum_{p\leq x} E(x/p) \log p \  = \ \text{Error} \tag{23}
$$
where  the Error is bounded by a multiple of $x$.  Dividing through by $x\log x$ we obtain
$$
\frac{E(x)}{x} \ =  \ - \frac 1{\log x}\  \sum_{p\leq x} \ \frac{E(x/p)}{x/p}  \cdot \frac{\log p}p +\ \text{Error},
$$
where the Error $\to 0$ as $x\to \infty$.
It is a little difficult to appreciate what this tells us. The right hand side can be viewed as $-1$ times the suitably weighted average of $E(t)/t$ for $t\leq x/2$ (use (7) to see that this really is a weighted average).  But this says that $E(x)/x$ is minus the average of $E(t)/t$, which is only consistent if that average is $0$, and therefore, we would hope to deduce that  $E(x)/x\to 0$ as $x\to \infty$, as desired.  One can make this deduction if one can prove that $E(x)/x$ does not change value quickly as $x$ varies, which is not straightforward. On the other hand, this argument is easily adapted to deduce that
$$
\liminf_{x\to \infty} \ \frac{E(x)}{x} \  =   \ - \  \limsup_{x\to \infty} \ \frac{E(x)}{x}  \ .
$$
This is as far as Selberg had gone using (19) when Erd\H os heard about Selberg's formula and started to work from it. Indeed both Erd\H os and Selberg went on   to deduce the prime number theorem using entirely elementary methods.\footnote{An unfortunate and unpleasant controversy arose as to who deserved credit for this first elementary proof of the prime number theorem. The establishment (as represented by the opinions of Weyl) judged that Erd\H os had ``muscled in'' on Selberg's breakthrough, that Selberg would have found the route to the elementary proof in time by himself. However Goldfeld [G4] provides an account of the controversy in which one cannot help but be sympathetic to Erd\H os. To my mind,
the controversy reflects two different perspectives on what is appropriate when one hears about the latest research of others, and what is not. Moreover what is appropriate changes over time and I do not think anyone would have questioned Erd\H os's behaviour today, nor would have been so unkind as Weyl. (See also [B1] and [Gr].)} We will not describe their proofs here since we will now  take these ideas in a different direction.

\subhead 9. Mean values of multiplicative functions \endsubhead We explained in section 2 that the prime number theorem is equivalent to the statement that the mean value of $\mu(n)$, the M\H obius function, for $n$ up to $N$, tends to $0$ as $N\to \infty$ (which is formulated in (5)). The beauty of reformulating the prime number theorem like this is that $\mu(n)$ is a multiplicative function, and this opens up many possibilities.  A {\sl multiplicative function} $f$ is one for which $f(mn)=f(m)f(n)$ whenever $m$ and $n$ are coprime integers.  Other important examples include

 $\cdot$ \ \ $n^{it}$, for fixed $t\in \Bbb R$;

 $\cdot$ \ \ $\chi(n)$, where $\chi$ is a Dirichlet character, which comes up when one studies arithmetic progressions;

 $\cdot$ \ \ $\tau(n)$, which counts the number of divisors of $n$;

 $\cdot$ \ \ $\sigma(n)$, the sum-of-divisors function, which arises when studying perfect numbers;

\noindent etc.  In all these cases we might ask for the function's mean value as we take the average up to infinity; that is
 $$
\lim_{N\to \infty} \ \frac 1N \ \sum_{n\leq N} f(n) .
$$
One should ask first whether this limit exists  and, if so, whether we can   determine its value. And, more importantly for the prime number theorem, can we come up with a simple classification of those multiplicative functions which have mean value $0$?

In fact the mean value, up to $N$, of $\chi(n)$  tends to $0$ as $N\to \infty$, of $\tau(n)$  tends to $\log N$, and of $\sigma(n)$ tends to $cN$ for some non-zero constant $c$. The most  interesting of our examples is $n^{it}$ with $t\ne 0$, since its mean value is
$$\frac 1N  \sum_{n=1}^N n^{it} \ \approx \   \frac 1N  \int_0^N u^{it} dt\ = \ \frac{N^{it}}{1+it}.$$
That is, the mean value does not tend to a limit as $N\to \infty$, but rather rotates steadily  around a circle of radius $1/\sqrt{1+t^2}$. We see here that the period works on a logarithmic scale, that is, we get roughly the same mean value for $N$ and $Ne^{2\pi/|t|}$.

In 1971 Hal\'asz resolved the key issue of determining which multiplicative functions do {\sl not} have mean value tending to $0$.  Restricting attention to multiplicative functions $f$ for which   $|f(n)|\leq 1$ for all $n$, there is the obvious example $1$, or any example much like $1$ (e.g. in which we perturb the value at each prime by just a small amount).  There is the generalization $n^{it}$,  for any real number $t$ (as we have just seen), and any small perturbations of that. Hal\'asz  proved  that these are essentially all the examples: The only multiplicative functions whose mean value do not tend to $0$ are ones that look a lot like $n^{it}$ for some real number $t$, that is, {\sl pretend} to be $n^{it}$. His proof  involves Dirichlet series to the right of $1$ and Parseval's identity, but never uses analytic continuation.

We now apply this to (5). If $\mu(n)$ does not have mean value $0$, then Hal\'asz's theorem tells us that $\mu(n)$ must {\sl pretend} to be $n^{it}$  for some real number $t$. Hence $\mu(n)^2$ must pretend to be $n^{2it}$, which implies that $t=0$, and hence $\mu(n)$ pretends to be $1$, a contradiction.  Formulating what ``pretends'' means takes a little bit of doing:\  If $f$ and $g$ are multiplicative functions with absolute value $1$, then $f$ pretends to be $g$ (meaning that they are not too different, at least in an appropriate average sense) if and only if  $f\overline{g}$ pretends to be $1$. Since the values of a multiplicative function only depend  on its values at primes and prime powers, we can restrict our attention to these. Now if $h$ pretends to be $1$ then we might measure that, by seeing how small $|1-h(p)|$ is, averaged in some way over the primes, or even $1-\text{Re }h(p)$, which turns out to be more natural (because of (31), below). Thus we define the {\sl distance}, ${\Bbb D}(f,g;x)$, between $f$ and $g$ for $n\leq x$, by:
$$
{\Bbb D}(f,g;x)^2 = \sum_{p\le x} \frac{1-\text{Re }(f(p)\overline{g(p)})}{p}.
$$
We say that $f$ is {\sl $g$-pretentious} if ${\Bbb D}(f,g;x)$ is less than some small constant. This allows us to formalize what we meant earlier (e.g. after (17)) when we wrote ``$p^{it}\approx -1$ for most primes $p$'' --- now we simply write that ${\Bbb D}(\mu, n^{it};x)$ is bounded.

${\Bbb D}(f,g;x)$ is not truly a distance; for example, it is only $0$ if  $f=g$ and, {\sl also}, $|f(n)|=1$ for all $n$. But for us it is more important that our notion of distance satisfies the triangle inequality
$$
{\Bbb D}(f,h;x) \ \leq \ {\Bbb D}(f,g;x) + {\Bbb D}(g,h;x). \tag{29}
$$
In particular we deduce that ${\Bbb D}(f^2,1;x)\le 2\ {\Bbb D}(f,\mu;x)$, and from this we easily deduce the prime number theorem. We are still looking for an elegant proof of the triangle inequality -- more on that in Appendix 2.

There is a direct connection between $\Bbb D(f,n^{it};x)$ and the Dirichlet series
$F(s):=\sum_{n\geq 1} f(n)/n^s$:\ If $\sigma=1+\frac 1{\log x}$ then
$$
|F(\sigma+it)| \ \asymp \  \log x \ \exp\left( - \Bbb D(f,n^{it};x)^2 \right) . \tag{31}
$$
where the symbol ``$\asymp$'' means that the ratio of the two sides is bounded, above and below, by positive constants. Hal\'asz's Theorem gives an upper bound for the mean value of $f$ in terms of the minimum of $\Bbb D(f,n^{it};x)$ as we range over $t$ in some box, $|t|\leq T$, where $T$ is a power of $\log x$ (that is, the minimum occurs at that $t$, with $|t|\leq T$, for which $n^{it}$ is ``closest'' to $f(n)$). From (31) this $t$ can also be thought of as the value at which $|F(\sigma+it)|$ is largest (up to a constant) out of  those $t$ for which $|t|\leq T$.

Besides the distance function, another key tool in working with mean values of multiplicative functions is a generalization of the argument we used in obtaining (7).
Now we are interested in $S(x):=\sum_{n\leq x} f(n)$. The trick is to evaluate
$$
\sum_{n\leq x} f(n) \log n \tag{37}
$$
in two ways, one analytic, the other algebraic. First analytically, note that (37) equals
$$
S(x) \log x - \int_1^x \frac{S(t)}t dt .
$$
Since $|S(t)|\leq \sum_{n\leq t} |f(n)|\leq t$, the integral here is $\leq \int_1^x 1 dt= x$, so the value is $S(x) \log x$ plus an error bounded by $x$, in absolute value. The second way to evaluate (37)  involves again writing $\log n$ as the sum of the logarithms of its prime and prime power divisors. As before we can bound the contribution of the prime power divisors, and so we are left with the ``identity''
$$
S(x) \log x \ = \  \sum_{p\leq x} f(p) \log p \ S(x/p) \ + \ \text{Error} \tag{41}
$$
where  the Error is bounded by a multiple of $x$. If we look at the special case $f=\mu$ and write $M(x)=\sum_{n\leq x} \mu(n)$ then we obtain
$$
M(x) \log x \ + \  \sum_{p\leq x} M(x/p) \log p \  = \ \text{Error} 
$$
where, analogously,  the Error is bounded by a multiple of $x$. Notice that this is {\sl exactly} the same functional equation as (23), with $M(.)$ replaced by $E(.)$ (which was the error term in the prime number theorem). This was a lot easier to derive than (23), and from here we can also prove that $M(x)/x\to 0$ as $x\to \infty$, much like Erd\H os and Selberg deduced that $E(x)/x\to 0$ as $x\to \infty$ from (23). This yields another ``elementary'' proof of the prime number theorem.

\subhead 10. What else do we count about primes? \endsubhead If one writes down the primes, it soon appears as if there are roughly equal numbers that end in a $1, 3, 7$ or $9$; in other words, in each residue class $a\pmod {10}$ where $(a,10)=1$. And it appears that, in general, for each fixed $q$, there should be roughly equal numbers of primes in each arithmetic progression $a \pmod q$ for any positive integer $a$ with $(a,q)=1$. However, even proving that there are infinitely many primes in any such arithmetic progression is a rather tough challenge. It was only in 1837 that Dirichlet did so, showing that the primes are equidistributed\footnote{By ``equidistributed'' we mean that there are roughly the same number of primes, up to $x$, in each arithmetic progression $a\pmod q$ with $(a,q)=1$.} in the arithmetic progressions if one weights them with a  $1/p$. In doing this, Dirichlet invented a generalization of the Riemann-zeta function, called a {\sl Dirichlet $L$-function}:\footnote{The astute reader might ask how Dirichlet could ``generalize'' the  Riemann-zeta function, 22 years before Riemann's paper! The fact is that $\zeta(s)$ was considered at length by Euler about one hundred years before Dirichlet; it was later named after Riemann, in honour of his trailblazing work.}  {\sl Dirichlet characters}, $\chi \pmod q$,  are multiplicative functions $\chi:\Bbb Z\to \Bbb C$ that are periodic with  minimal period $q$, and $\chi(n)=0$ if $(n,q)>1$. The most interesting are the real characters (i.e.\ those characters that can only take the values $-1,1$ and $0$, and do, in fact, take each of those values) like the Legendre symbol $\left( \frac{.}q\right)$. For each Dirichlet character we create the {\sl Dirichlet $L$-function}
$$
L(s,\chi) \ := \ \sum_{n\geq 1} \ \frac{\chi(n)}{n^s},
$$
which is absolutely convergent when Re$(s)>1$. This can be analytically continued to the whole complex plane with no poles.  Dirichlet's proof is an elegant piece of combinatorics which easily leads to his theorem that there are infinitely many primes $\equiv a \pmod q$ whenever $(a,q)=1$ provided one can prove that
$$
L(1,\chi)\ne 0 \ \ \text{for all real characters} \ \chi \pmod q.
$$
This is still a whole lot harder to prove than one might guess, even though there are many different proofs. The most interesting proof, due to Dirichlet himself, shows that $L(1,\chi)$ can be determined as a simple multiple of the size of a certain group that comes up in algebra. This {\sl Dirichlet's class number formula} was the first deep connection found between algebra and analysis and is the pre-cursor of so many of the great theorems and conjectures of the last thirty years in number theory.\footnote{Like Wiles' Theorem, the Birch-Swinnerton Dyer conjecture, etc.}

The proof of the prime number theorem was soon modified, using Dirichlet series,  to show that,  whenever $(a,q)=1$ and $(b,q)=1$,
$$
\lim_{x\to\infty} \ \frac{ \#\{ p\leq x:\ p \ \text{prime, and}\ p\equiv a \pmod q\} }
{ \#\{ p\leq x:\ p \ \text{prime, and}\ p\equiv b \pmod q\} } \ \ \text{exists and equals} \ 1;
$$
that is, the primes are equidistributed amongst the plausible arithmetic progressions mod $q$. This is called the {\sl Prime Number Theorem for arithmetic progressions mod} $q$.

Approaching this with the M\H obius function, one can show that the prime number theorem for arithmetic progressions holds mod $q$  if and only if
$$
\lim_{N\to \infty} \ \frac 1N \ \sum_{n\leq N} \mu(n)\chi(n) \ = \ 0
$$
for every Dirichlet character $\chi \pmod q$. As in the classical proof, this is easily proved using Hal\'asz's Theorem when $\chi$ takes complex values (since then $\chi(p)\mu(p)$ is often quite far from the real line). If  $\chi$ is a real character it also follows immediately from Hal\'asz's Theorem provided   $L(1,\chi)\ne 0$. So the ``pretentious proof'' hinges on the same issue as the classical proof.

\subhead 11. Is a pretentious proof, an elementary proof? \endsubhead
The prime number theorem can be phrased as: \ For all $\epsilon>0$, we have $|\theta(x)-x| < \epsilon x$ once $x$ is sufficiently large.  However we believe that much more is true, namely $|\theta(x)-x| < x^{1/2}(\log x)^2$, so the question becomes how close we can get to our belief.  We have seen that improving the error term in the prime number theorem is equivalent to exhibiting wider regions to the left of the ``$1$-line'' that contain no zeros of $\zeta(s)$, so-called {\sl zero-free regions}. Proving such results requires complicated analysis of various explicit formulas involving the zeros of $\zeta(s)$  
(details can be found, for instance, in [Da] and [Bo]).  The key idea is that we understand  the values of $\zeta(s)$ well to the right of the $1$-line, and we can use that understanding, via such explicit formulas, to get some control over $\zeta(s)$ just to the left of the $1$-line (we write ``just'', since the zero-free regions that have been obtained are so very narrow). These are beautiful and subtle proofs but they give relatively weak results. Moreover, the main application is to discuss issues about prime numbers that are, essentially, questions that arise in the elementary world to the right of the $1$-line. It seems strange to work so hard to extrapolate our knowledge of $\zeta(s)$ to the right of the $1$-line, in order to get a meagre understanding just to the left of the $1$-line, so as to answer questions to the right of the $1$-line. Why on earth should we cross the $1$-line at all? The goal of these methods is to render that journey unnecessary. 

Pretentious methods use non-trivial techniques of complex analysis (in particular Perron's formula, as we will see), but not analytic continuation, nor Cauchy's Theorem and residue computations, nor subtle calculations of zeros of analytic continuations. On the other hand, the calculations used in pretentious techniques can be challenging, though usually they can be reduced to completely elementary techniques, at the cost of further complications. So, technically, one could invoke Selberg's definition to say that these are elementary techniques, though that misses the point. The main issue is that these methods avoid needing the Dirichlet series $F(s)$ to be analytically continuable,\footnote{Which is a rare property for a Dirichlet series, a technical property that sometimes seems unreasonably convenient, though it  is conjectured to be true for most $L$-functions of direct arithmetic interest.} and so are much more widely applicable.

\subhead 12. Primes in arithmetic progressions (questions involving uniformity) \endsubhead
If we begin computing primes in arithmetic progressions mod, say, $101$, we notice that, quite soon, the primes are roughly equidistributed in all of the $100$ possible progressions.\footnote{For example, by the time there are $100$ primes, on average, in each arithmetic progression mod $101$,  the least is $87$ primes in some arithmetic progression mod $101$, and the most is $109$.  By the time there are $1000$ primes, on average, in each arithmetic progression,  the least is $968$ and the most is $1030$. By the time there are $10,000$ primes, on average, in each arithmetic progression,  the least is $9912$ and the most is $10070$.}
So there is an important, new question for primes in arithmetic progressions:
\smallskip

\centerline{\sl When can we expect roughly equal numbers of }

 \centerline{\sl  primes in each arithmetic progression mod $q$?}
\smallskip

\noindent After a lot of  computing, researchers guess that there should be roughly equal numbers of primes up to $x$, in each arithmetic progression $a\pmod q$, for each $a$ with $(a,q)=1$, once $x>c_\epsilon q^{1+\epsilon}$, for some constant $c_\epsilon$, for any fixed $\epsilon>0$. This is far out of reach of what we can prove. Indeed,  the best result we have a plan of how to prove is  that the primes are roughly equidistributed mod $q$ once  $x>c_\epsilon q^{2+\epsilon}$. This plan, though, involves proving the {\sl Generalized Riemann Hypothesis},\footnote{That is, if $\rho$ is a zero of any Dirichlet $L$-function with $0\leq \text{Re}(\rho)\leq 1$, then Re$(\rho)=\frac 12$.} which seems very far out of reach.

So what can we prove unconditionally? In both the classical theory, and the pretentious approach, the issue is how close $L(1,\chi)$ is to $0$ for real characters $\chi$. A lower bound  of the shape $L(1,\chi)>c/\sqrt{q}$ follows immediately from Dirichlet's class number formula, and  leads to the result that the primes are roughly equidistributed mod $q$ once  $x>c e^{\sqrt{q}}$ for some constant $c>0$ that one can determine. This lower bound for $x$ is far bigger than what we expect to be true.

In 1936 Siegel improved this lower bound to: For all $\epsilon>0$ there exists a constant $c_\epsilon>0$ such that $L(1,\chi)>c_\epsilon/q^\epsilon$, and so
the primes are roughly equidistributed mod $q$ once  $x>\kappa_\epsilon e^{q^\epsilon}$.
But there is a catch. The method of proof does not allow one to determine $c_\epsilon$: Note that we   are not saying  that it has not been computed, but rather that   {\sl it cannot be computed}. The proof is very surprising in that Siegel splits his considerations into two complementary cases:
\smallskip
Either, the Generalized Riemann Hypothesis is {\sl true}, so it is easy to compute $c_\epsilon$.

Or, the Generalized Riemann Hypothesis is {\sl false}, in which case it is easy to compute $c_\epsilon$ in terms of any given counterexample.
\smallskip

\noindent The problem with this dichotomy is the second case. The Generalized Riemann Hypothesis is unresolved, and were the Generalized Riemann Hypothesis to be false,  then Siegel's proof can only provide a constant  once some counterexample to the Generalized Riemann Hypothesis is known.

All this talk of Riemann Hypotheses in the proof of Siegel's theorem means that we are involving zeta-functions to the left of the $1$-line, and so I had believed that this result could only be obtained by classical means. That was my prejudice, until my postdoc, Dimitris Koukoulopoulos,\footnote{And now my faculty colleague.}  came up with a very subtle elementary argument that allowed him to completely replace Siegel's argument by a purely pretentious one, with no analytic continuations in sight. One can find links between his proof and that of Siegel's (as developed by Pintz [P1]) and so, rather amazingly, we now have an ``elementary proof'' of Siegel's theorem.

Koukoulopoulos [K1] also showed that  the primes are not only equidistributed mod $q$ once  $x>\kappa_\epsilon e^{q^\epsilon}$, but that the ratio is {\sl very} close to $1$ (as in the {\sl Siegel-Walfisz Theorem}). This in turn allows one to use the large sieve to prove the {\sl Bombieri-Vinogradov Theorem}. The Bombieri-Vinogradov Theorem can be interpreted as stating that the primes are more-or-less equidistributed mod $q$ {\sl for almost all} $q<x^{1/2-\epsilon}$. This is the consequence we expect from the Generalized Riemann Hypothesis, but we obtain it only for most $q$, not necessarily all $q$.

\subhead 13. Pretentiousness is repulsive\endsubhead Seemingly, one of the deepest results about $L$-functions is that their zeros ``repel'' each other. That is, they do not like to be too close together. In particular one cannot have zeros of two Dirichlet $L$-function both close to $1$, and this can be re-phrased as saying that there is at most one real character $\chi \pmod q$, amongst all the real characters with $q$ in the range $Q<q\leq 2Q$, for which $L(1,\chi)< c/\log Q$. Hence $L(1,\chi)\geq c/\log q$ for all of the other real characters with $q$ in this range, and therefore one can state a strong prime number theorem for the arithmetic progressions for all these other moduli.\footnote{I have simplified here a little bit, rather than get in to the technicalities of primitive and induced characters. For more on this, see Davenport's book [Da].}  In fact with such a strong lower bound on $L(1,\chi)$ one can show that the primes are roughly equidistributed mod $q$ once  $x>q^A$ for some sufficiently large $A$ (how large depends on the constant $c$, and how nearly you want the primes to be equidistributed).

Rather surprisingly,   these repulsion results are much easier to prove in the pretentious world. Basically $L(1,\chi)$ being very small means that $\mu$ is $\chi(n) n^{it}$-pretentious for some real number $t$, and so if $L(1,\psi)$ is also small then $\mu$ is $\psi(n) n^{iu}$-pretentious for some real number $u$. Now if $\mu$ is very close to $\chi(n) n^{it}$ as well as to $\psi(n) n^{iu}$, then they are close to each other (which formally follows from our  triangle inequality), and therefore the Dirichlet character
$\overline{\psi}\chi$ is $n^{i(u-t)}$-pretentious, which is easily shown to be impossible. This all goes to show that {\sl pretentiousness is repulsive}.

\subhead 14. The pretentious large sieve \endsubhead Perhaps the deepest proofs in the classical analytic number theory approach to the distribution of prime numbers, are the proofs of {\sl Linnik's Theorem}; that is that there exist constants $c>0$ and $L>0$ such that for any positive integers $a$ and $q$, there is a prime $\leq cq^L$ which is $\equiv a \pmod q$. Linnik's 1944 proof [L1] has been improved many times (e.g.~in Bombieri's [Bo]) but remains delicate and subtle.  Inspired by a new, technically elementary proof in November 2009 given by Friedlander  and Iwaniec [F4], Soundararajan and I went on  to develop an idea we had for a {\sl pretentious large sieve} [GS1], and we ended up giving what is surely the shortest and technically easiest proof of Linnik's Theorem, though bearing much in common with an earlier proof of   Elliott [E2].   This new technique has enormous potential, because it can replace some very subtle classical techniques, and yet does not require the function involved to have an $L$-function that can be analytically continued. We made one application, with de la Br\'eteche [GS2], to better  understand the solutions to Pythagoras' equation $a^2+b^2=c^2$ mod $p$ (as well as to several other additive number theory problems). With Adam Harper we can now prove a weak form of Hoheisel's deep theorem on primes in short intervals: That is, there exists a constant $\delta>0$ and a constant $c_\delta>0$ such that if $x^{1-\delta}<y\leq x$ then 
$$
\#\{ p \ \text{prime}: \ x<p\leq x+y\}  \ \geq \ c_\delta\ \frac{y}{\log x}.
$$

\subhead 15. From a collection of {\sl ad hoc} results, to a new approach to prime numbers \endsubhead
Our easy proof of Linnik's Theorem suggested to Soundararajan and me that we should be able to prove all of the basic results of analytic number theory {\sl without ever using analytic continuation}. Since early 2010 we have been working on developing this new approach. Our goal is to reprove all of the key results in the standard classical books [Da] and [Bo] using only ``pretentious methods''.  Within a year we found that we could prove some version of all of the results, perhaps not as strong, but much the same in principle. In doing this we have stood firmly on the shoulders of giants. Many analytic number theorists have developed ideas about multiplicative functions, over the last 50 years, that have allowed them to prove results on different aspects of prime numbers. Those who have been most central to our description of the subject are Erd\H os, Selberg, Wirsing, Bombieri, Delange, Daboussi, Hildebrand, Maier, Hall, Pomerance,Tenenbaum, Pintz, Elliott, Montgomery and Vaughan,  Friedlander, Iwaniec and Kowalski,\dots

Despite being able to prove some version of all of the principal results known on the distribution of primes, we increasingly found ourselves frustrated for three reasons:

\medskip

\item{(1)}\ Although we could show the prime number theorem, we could not show that convergence is anywhere like as fast as had been shown by classical means. In fact we {\sl could not} see how one might use pretentious methods to even prove something as (relatively) weak as $|\theta(x)-x|\leq x/\log x$ for $x$ sufficiently large.

\item{(2)}\ There are many strong results in the subject that are proved assuming the Riemann Hypothesis. We could not conceive of proving analogous results since the  Riemann Hypothesis is a conjecture about the zeros of the analytic continuation of $\zeta(s)$, something we are trying to avoid discussing at all.

\item{(3)} All proofs of Hal\'asz's Theorem, which lies at the center of the whole theory, were not only complicated  but  also hard to motivate. We had modified this proof in several ways, for example in the proof of the pretentious large sieve, and this led to a lot of the theory seeming somewhat obscure, even if technically straightforward.

\medskip

\noindent Fortunately several of the best young people in analytic number theory got interested in these issues, and they have satisfactorily resolved all of them, as I will now describe:

\subhead 16. The strongest known form of the prime number theorem \endsubhead The prime number theorem can be phrased as $\theta(x)/x\to 1$ as $x\to \infty$. The proofs of 1896 immediately yielded that for any fixed $A>0$ we have
$$
|\theta(x)-x| \leq \frac x{(\log x)^A}
$$
for all sufficiently large $x$.  In fact de la Vall\'ee Poussin proved the much stronger result that
there exists a constant $c>0$ for which
$$
|\theta(x)-x| \leq x / \exp\left({c\sqrt{\log x}}\right)
$$
for all sufficiently large $x$.  The strongest version proved unconditionally is from 1959 (by Korobov and Vinogradov), and gives
$$
|\theta(x)-x| \leq x / \exp\left( c \frac{ (\log x)^{3/5} }  { (\log\log x)^{1/5} } \right) .  \tag{43}
$$
There has been no improvement on this in over 50 years, yet it is so far from what we believe to be true, and can prove  assuming the Riemann Hypothesis.

In directly using Hal\'asz's Theorem, applied to the multiplicative function $\mu(n)$ as described above, one can prove results like
$|\theta(x)-x| \leq x / (\log x)^\tau$ for sufficiently large $x$, for some $\tau<1$, and one can show that it is impossible to obtain  larger $\tau$ as an immediate corollary. This is a lot weaker result than the simplest results that one obtains from classical methods.

Selberg showed that if $f(p)=\alpha$ for all primes $p$ (with $|\alpha|\leq 1$), then the bounds given by Hal\'asz's Theorem are very near to the truth, unless $\alpha=0$ or $-1$ in which case the mean value is much smaller.  This includes the example of the mean value of  $\mu(n)$. In fact, Delange went on to show that if the mean value of $f(p)$ is $\alpha$ then much the same result holds; and hence the mean value of $\mu(n)\chi(n)$ converges rapidly to $0$ for any Dirichlet character $\chi$. The {\sl Koukoulopoulos converse theorem} [K2] goes one big step forward, stating  that if the mean value of $f(n)$ is small then  $f(p)$ must average
$0$ or $-1$ over the primes. This opens the door to getting much stronger upper bounds for the mean value of $f$, via appropriate modifications of  Hal\'asz's Theorem. Indeed Koukoulopoulos was then able to prove the strongest known versions of both the prime number theorem, as in (43), and the prime number theorem for arithmetic progressions, using only pretentious methods, never venturing to the left of  the $1$-line.

At first sight it is surprising that he could not do better than the classical proofs. After all, if   Koukoulopoulos's proofs are so different from the classical proofs, then why would he also come up with such an unlikely bound as in (43)?  The reason is that, despite appearances, the proofs are fundamentally the same. The classical proof uses deep tools of analysis, 
which are stripped away in Koukoulopoulos' proof, suggesting that the use of zeros of $\zeta(s)$ is artificial.

\subhead 17. The Pretentious Riemann Hypothesis\endsubhead The Riemann Hypothesis tells us that the zeros of (the analytic continuation of) $\zeta(s)$ are far into the domain of analytic continuation, that is, they are all on the ``$\frac 12$-line''.  Does $\zeta(s)$ feel the effect of this to the right of the ``$1$-line''? Can we recognize the Riemann Hypothesis to the right of the $1$-line?

To count primes we looked at $\zeta'(s)/\zeta(s)$. The Riemann Hypothesis is equivalent to this not having any poles, other than at $s=1$, to the right of the $\frac 12$-line. One can remove the pole at $s=1$  by working instead with 
$$\frac{\zeta'(s)}{\zeta(s)} + \frac 1{s-1} \ ,$$ 
or even
 $$\frac{\zeta'(s)}{\zeta(s)} + \zeta(s) \ . \tag{47}$$
If this function's Taylor series converges around any given point $s_0$ to the right of the $1$-line,  within the ball 
$$ B\left(s_0,\frac 12\right):= \left\{ s:\ |s-s_0| < \frac 12 \right\} ,$$
 then  there can be no poles in that region. The union of all those balls equals the domain to the right of the $\frac 12$-line; that is
$$
\bigcup_{s_0: \ \text{Re}(s_0)>1}  B\left(s_0,\frac 12\right) \ = \ \left\{ s:\ \text{Re}(s)>\frac 12 \right\} .
$$
Therefore, we have proved that the Riemann Hypothesis holds if the Taylor series for (47) at $s=s_0$ converges within $B\left(s_0,\frac 12\right)$, for any $s_0\in \Bbb C$ to the right of the $1$-line.

The $k$th coefficient of the Taylor series for $f(s)$ at $s=s_0$ is given by $f^{(k)}(s_0)/k!$. We can therefore guarantee that the Taylor series   converges absolutely within $B\left(s_0,\frac 12\right)$, if  $|f^{(k)}(s_0)|\leq c(s_0) k!2^k $ for every integer $k\geq 0$, for some constant $c(s_0)$ which may depend on $f$ and $s_0$. We conjecture such a hypothesis for the function in (47):

 \proclaim{The Pretentious Riemann Hypothesis}\   For all $\epsilon>0$ there exists a constant $c_\epsilon>0$ such that for every integer $k\geq 1$ we have
$$
\Big| \Big( \frac{\zeta'}{\zeta} (s)+\frac{1}{s-1}\Big)^{(k)} \Big| \leq c_\epsilon k! 2^k (1+t^\epsilon)
$$
uniformly for $s=\sigma+it$ with $1\leq \sigma<2$ and $0\leq t\leq e^k$.
\endproclaim

Using Koukoulopoulos' methods one can show that if the pretentious Riemann Hypothesis is true then  $|\psi(x)-x|<\kappa_\epsilon x^{1/2+\epsilon}$ for any fixed $\epsilon>0$, which in turn implies the Riemann Hypothesis. On the other hand, the  Riemann Hypothesis implies\
$$
\Big| \Big( \frac{\zeta'}{\zeta} (s)+\zeta(s) \Big)^{(k)} \Big| \leq c \ k! 2^k \log t
$$
for $s=\sigma+it$ with $\sigma\geq 1$ and  $t\geq 1$, which is somewhat more than we asked for in the pretentious Riemann Hypothesis. Together these remarks imply:
\bigskip

\centerline{\sl The Riemann Hypothesis holds if and only if the Pretentious Riemann Hypothesis holds.}

\subhead 18. A re-appraisal of the use of Perron's formula \endsubhead Soundararajan and I had written up as palatable a  proof as we could  of Hal\'asz's Theorem for the first drafts of our book [GS1], but even we had to admit that it was difficult to motivate.  So I was delighted when, in January 2013, my new postdoc, Adam Harper,  suggested a new, simpler path to a proof of Hal\'asz's Theorem. Subsequently we have developed his idea with him, and find ourselves re-appraising the use of Perron's formula when summing coefficients of Dirichlet series.

In Riemann's approach one takes the formula (11), and shifts the line of integration far to the left side of the complex plane. In the pretentious approach one stays with the same line of integration. But then how can one get an accurate estimate, or even a decent upper bound, since the integral of the absolute value of the integrand is usually much larger than the value of the integral? There are several important observations involved.  First though, let's look at this in more generality, with the identity
$$
\sum_{n\leq x} f(n) \ =  \ \frac 1{2i\pi} \int_{\text{Re}(s)=\sigma} F(s) \frac{x^s}s ds, \tag{53}
$$
where $F(s):=\sum_{n\geq 1} f(n)/n^s$ and $\sigma=1+\frac 1{\log x}$.   One can give a version of this (like in the proof of the prime number theorem), with the values of $s$ running over only those $t$ where $|t|\leq T$, for some suitably chosen $T$ (taking $T$ as a power of $\log x$ will do). Then we can take absolute values in the integrand, noting that $|x^s|=x^\sigma=ex$ to get an upper bound
$$\align
\frac 1{2\pi}  \int\Sb s=\sigma+it \\ |t|\leq T \endSb   |F(s)| \frac{|x^s|}{|s|}  dt \ &\leq\ 3x \cdot \max_{|t|\leq T} |F(\sigma+it)| \cdot  \frac 1{2\pi} \int_{|t|\leq T} \frac{1}{|\sigma+it|} dt \\
&\leq \  x \cdot \max_{|t|\leq T} |F(\sigma+it)| \cdot   (\log T+1) .
\endalign
$$
This   would more-or-less be Hal\'asz's Theorem (via (31)) if the upper bound was divided through by $\log x$. This is encouraging since our approach in getting this upper bound was very crude, and we can surely refine it a bit.

Studying the integrand $F(s)x^s/s$, we might expect that $F(s)x^\sigma/s$ does not change much while $x^{it}$ rotates once around the unit circle (which requires an interval, for $t$, of length  $2\pi/\log x$).  The easiest way to pick up this cancellation is to integrate by parts, so that (53) becomes:
$$
\ =  \ \frac 1{2i\pi} \int_{\text{Re}(s)=\sigma} F(s) \frac{x^s}{s^2 \log x}
\ - \  \frac 1{2i\pi} \int_{\text{Re}(s)=\sigma} F'(s) \frac{x^s}{s \log x} .
$$
The ``$\log x$'' in the denominator is the cancelation. The first term, because of the  $s^2$ in the denominator, is sufficiently small to be ignored, and so we are left with
$$
-    \frac 1{2i\pi \log x} \int_{\text{Re}(s)=\sigma} F'(s) \frac{x^s}{s} \ = \
-    \frac 1{2i\pi \log x} \int_{\text{Re}(s)=\sigma} \frac{F'(s)}{F(s)} \cdot F(s) \frac{x^s}{s}.  \tag{59}
$$
If we take absolute values here, much as we did in (53), then we get the desired bound so long as $|F'(s)/F(s)|$ is ``small'' in a certain average sense. It is indeed this small for many multiplicative functions $f$ of interest to us, but not all, so another idea is needed.

This is where the key new idea of Harper comes in.
Going from (53) to (59), we gained a factor of $\frac 1{\log x} \cdot \frac{F'(s)}{F(s)}$, which does improve things and  gaining another such factor would be enough to get us to our goal. We cannot quite do this, but a variant gets us there, using the flexibility of slightly varying the (vertical) line of integration:
$$
\int_0^\Delta   \left( \frac{1}{i\pi} \int_{1-iT}^{1+iT}
 \frac{F'}{F}(s-\beta) \cdot \frac{F'}{F}(s+\beta) \cdot
F(s+\beta) \frac{ x^{s-\beta}}{s+\beta} ds  \right) d\beta
$$
where $\Delta$ is a suitably chosen multiple of $1/\log x$.
There are similarities between this and   the  formula  used in the  usual proof of Hal\'asz's Theorem, but it is now much clearer how we got  here (which means that this new technique is much more flexible). In particular it allows us to obtain asymptotic formulae if the mean value is ``reasonably well-behaved.''

\subhead 19. New results \endsubhead  Our goal in this project is to reprove all of the results of classical analytic number theory from our new perspective. This is a worthwhile project but typically mathematical researchers look forward to what comes next, not to what has been. So to truly justify developing these methods, one might ask whether we can prove results that classical methods could not?

In fact, these pretentious techniques were not born from trying to reprove old results, but rather from proving new results on old problems and seeing a pattern emerge in our proofs. The inspiration was the first ``big'' improvements in bounds for sums of characters in ninety years [GS3]. By understanding that a character sum could be large only if the character pretended to be a different character with much smaller modulus, we were able to find our improvement. Moreover we found new inter-relations between large character sums that had not previously been known to exist, or even guessed at.

There have been other results: other questions on character sums [GS4, Go], large $L$-function values [GS5], least non-residues [XL], convexity problems for $L$-functions [S1], and most spectacularly, Soundararajan's work, with Holowinsky, completing the proof of Arithmetic Quantum Unique Ergodicity [So1, So2, HS], which had been a famous conjecture. Very recently, Matom\" aki and Radiziwill [MR] have used such techniques to study  sign changes in the coefficients of holomorphic Hecke cusp forms.

Hal\'asz's theorem is bound to be a better tool to study more general analytic problems than classical analytic methods since the Dirichlet series arising from the given multiplicative function do not need to be analytic, which was the whole point of using zeta-functions.

On the other hand zeta-functions have a rich history, and are central to many key themes in mathematics.  Dirichlet $L$-functions are the zeta-functions of weight one, the simplest class. Next come the weight two zeta-functions, which include the $L$-functions  associated to elliptic curves. There is much to do to establish a pretentious theory here. The classical theory can prove much less with these $L$-functions, so we can hope that pretentious techniques might have significant impact on arithmetic questions associated to these $L$-functions.

\subhead{Acknowledgements} \endsubhead 
Thanks are due to Eric Naslund and Dimitris Koukoulopoulos, as well as the anonymous referee and  Frank Farris the editor, who read and commented helpfully on a preliminary version of this article.

\subhead{Appendix One: Factorization Tables}  \endsubhead  In Gauss's December 24th, 1849 letter to Encke, he wrote:

\block {\eightpoint
    {\sl The 1811 appearance of Chernac's {\rm Cribrum Arithmeticum}    gave me great joy, and (since I did not have the patience for a continuous count of the series) I have very often employed a spare
unoccupied quarter of an hour in order to count up a chiliad here and there;
however, I eventually dropped it completely, without having quite completed
the first million.}
}
\endblock

\noindent Figure 1 is a photograph of Chernac's Table of Factorizations of all integers up to one million, published in 1811, which was used by Gauss. My colleague, Anatole Joffe, kindly presented his copy of these tables to me when he retired. Nowadays I also happily distract myself from boring office-work, by flicking through to discover obscure factorizations!

There are exactly one thousand pages of factorization tables in the book, each giving the factorizations of one thousand numbers. For example, page 677, seen here, enables us to factor all numbers between $676000$ and $676999$. The page is split into 5 columns, each in two parts, the ten half columns on the page each representing those integers that are not divisible by 2, 3 or 5, in an interval  of length 100.\footnote{Chernac trusted that the reader could easily extract those factors for him- or herself.} On the left side of a column is a number like $567$, which represents, on this page, the number $676567$ to be factored. On the right side of the column we see $619\cdot 1093$ which gives the complete factorization of $676567$. On the other hand for  $589$, which represents  the prime number $676589$, the right column  simply contains ``------'', and hence that number is prime. And so it goes for all of the numbers in this range. It only takes a minute to get used to these protocols, and then the table becomes   very useful if you do not have appropriate factoring software at your disposal.

\subhead{Appendix Two. A proof and a challenge} \endsubhead The triangle inequality, (29), lies at the heart of this new theory. There are now several proofs, none of which are particularly elegant. The best was given by Eric Naslund, as part of an undergraduate research project in 2011, which I will reproduce now. To prove (29) it suffices to prove the simpler inequality
$$
\eta(w,y) \le \eta(w,z) +\eta(z,y). \tag{61}
$$
where $\eta(z,w)^2 := 1-\text{\rm Re}(z\overline{w})$, for any $w, y, z\in \{ z\in \Bbb C:\ |z|\leq 1\}$.

\demo{Proof}  (Eric Naslund) Let $r=|z|$. Then, since $|w\overline{z}|\leq r$ and $|z\overline{y}|\leq r$, we can write $w\overline{z}=r(a+bi),\ z\overline{y}=r(c+di)$ and $ w\overline{y}=(a+bi)(c+di)$ where $a+bi$ and $c+di$ lie in the unit disk.

Note that $1+ra, 1+rc\leq 2$, so that $(1+ra)(1+rc)\leq 4$, and hence
$$
2\sqrt{1-ra}\sqrt{1-rc}\geq\sqrt{1-r^{2}a^{2}}\sqrt{1-r^{2}c^{2}}\geq\sqrt{1-a^{2}}\sqrt{1-c^{2}}\geq bd.
$$
Now, if $a+c\geq0$ then $1-ra-rc+ac\geq1-a-c+ac=(1-a)(1-c)\geq0$; and if
$a+c\leq0$ then $1-ra-rc+ac\geq1+ac\geq0$. Either way $1-ra+1-rc\geq 1-ac$.
Adding this to the last displayed equation and taking square roots of both sides, we obtain
$$\eta(w,z)+\eta(z,y)=\sqrt{1-ra}+\sqrt{1-rc}\geq\sqrt{1-ac+bd}=\eta(w,y).$$
\enddemo

We would still like to see a ``proof from the book'' \footnote{The great Paul Erd\H os, used to say that the Supreme Being has a book of all of the best proofs, and just occasionally we are allowed to glimpse at a page. When you have such a proof, it is obvious that it is from ``the book''! When he was still alive, there was no greater compliment than Erd\H os remarking, as he occasionally did, ``I think that is from the book.''} of (61), a more natural and easy proof. I will leave this as a competition for our readers. Please email me your proof. The best one will appear in our book, with appropriate credit.

\enddocument

Davenport's {\sl Multiplicative Number Theory} and Bombieri's {\sl Le Grand Crible}
\enddocument